\RequirePackage{ifthen}
\newboolean{MPA}
\setboolean{MPA}{false}

\ifthenelse {\boolean{MPA}}
{
\documentclass[numbook]{svjour3} 
\smartqed
\usepackage[margin=1.3in]{geometry}
} {
\documentclass[11pt]{article}
\usepackage[margin=1in]{geometry}
}

\usepackage[title]{appendix}

\usepackage{graphicx}
\usepackage{amsmath}
\usepackage{amssymb}
\usepackage{hyperref}
\usepackage{color}

\usepackage{cite} 

\usepackage{mathrsfs} 

\usepackage{xspace} 

\usepackage{soul} 

\usepackage[makeroom]{cancel} 

\usepackage{algorithm}
\usepackage[noend]{algpseudocode}

\floatname{algorithm}{}

\usepackage{accents}

\usepackage{enumitem} 
\usepackage[capitalize,noabbrev]{cleveref}
\crefname{question}{Question}{Questions}
\crefname{step}{Step}{Steps}
\crefname{claim}{Claim}{Claims}
\crefname{problem}{Problem}{Problems}
\crefname{observation}{Observation}{Observations}
\crefname{fact}{Fact}{Facts}

\usepackage{thmtools}
\usepackage{thm-restate}

\newcommand{\R}{\mathbb R}

\newcommand{\Q}{\mathbb Q}

\newcommand{\ind}{h}

\DeclareMathOperator{\ptw}{ptw}
\DeclareMathOperator{\itw}{itw}
\DeclareMathOperator{\tw}{tw}
\DeclareMathOperator{\inc}{inc}
\DeclareMathOperator{\prim}{prim}

\let\st\relax
\DeclareMathOperator{\st}{s.t.}

\DeclareMathOperator{\poly}{poly}

\newcommand{\A}{\mathcal A}
\newcommand{\I}{\mathcal I}

\newcommand{\abs}[1]{\lvert#1\rvert}
\newcommand{\absL}[1]{\left\lvert#1\right\rvert}

\newcommand{\pare}[1]{\left(#1\right)}
\newcommand{\bra}[1]{\left\{#1\right\}}
\newcommand{\sbra}[1]{\left[#1\right]}

\newcommand{\transp}{\mathsf T}

\newcommand{\inner}[2]{\left\langle #1, #2 \right\rangle}

\usepackage{comment}
\specialcomment{knownproof}{\color{blue}}{\color{black}}
\excludecomment{knownproof}

\newcounter{step}

\ifthenelse {\boolean{MPA}}
{

\newenvironment{prf}[1][]
{\begin{proof}}
{\qed \end{proof}}

\newenvironment{prfc}[1][]
{\begin{proof}{\emph{#1 }}}
{\qed \end{proof}}

\newcounter{claim} 
\renewenvironment{claim}[1][]
{\refstepcounter{claim} \begin{trivlist} \item[] {\bf Claim~\theclaim}\space#1 \itshape}
{\end{trivlist}}

\journalname{Mathematical Programming A}

} {

\usepackage{amsthm}
\newtheorem{theorem}{Theorem}
\newtheorem{proposition}{Proposition}
\newtheorem{lemma}{Lemma}

\newtheorem{corollary}{Corollary}

\newtheorem{claim}{Claim}
\newtheorem{example}{Example}

\newenvironment{prf}[1][]
{\begin{proof}}
{\end{proof}}

\newenvironment{prfc}[1][]
{\begin{proof}[Proof #1]}
{\end{proof}}

}

\DeclareMathOperator{\obj}{obj}

\newtheorem{observation}{Observation}
\newtheorem{fact}{Fact}

\begin{document}


\title{Factorized binary polynomial optimization}

\ifthenelse {\boolean{MPA}}
{

\author{Alberto Del Pia}
\institute{Alberto~Del~Pia \at
              Department of Industrial and Systems Engineering 
              \& Wisconsin Institute for Discovery \\
              University of Wisconsin-Madison, Madison, WI, USA \\
              \email{delpia@wisc.edu}}
}
{
\author{Alberto Del Pia
\thanks{Department of Industrial and Systems Engineering \& Wisconsin Institute for Discovery,
             University of Wisconsin-Madison, Madison, WI, USA.
             E-mail: {\tt delpia@wisc.edu}.}}
}

\date{July 4, 2024}

\pagenumbering{gobble} 

\maketitle


\begin{abstract}
In \emph{binary polynomial optimization,} the goal is to find a binary point maximizing a given polynomial function.
In this paper, we propose a novel way of formulating this general optimization problem, which we call \emph{factorized binary polynomial optimization.}
In this formulation, we assume that the variables are partitioned into a fixed number of sets, and that the objective function is written as a sum of $r$ products of linear functions, each one involving only variables in one set of the partition.
Our main result is an algorithm that solves factorized binary polynomial optimization in strongly polynomial time, when $r$ is fixed.
This result provides a vast new class of tractable instances of binary polynomial optimization, and it even improves on the state-of-the-art for quadratic objective functions, both in terms of generality and running time.
We demonstrate the applicability of our result through the \emph{binary tensor factorization problem,} which arises in mining discrete patterns in data, and that contains as a special case the \emph{rank-$1$ Boolean tensor factorization problem.}
Our main result implies that these problems can be solved in strongly polynomial time, if the input tensor has fixed rank, and a rank factorization is given.
For the \emph{rank-$1$ Boolean matrix factorization problem,} we only require that  the input matrix has fixed rank.

\ifthenelse {\boolean{MPA}}
{
} {}
\end{abstract}

\ifthenelse {\boolean{MPA}}
{}{
}

\section{Introduction}
\label{sec intro}

Binary polynomial optimization, i.e., the problem of finding a binary point maximizing a given polynomial function, is a fundamental problem in discrete optimization with a wide range of applications across science and engineering.
To formulate this optimization problem, we use the hypergraph representation from \cite{dPKha17MOR}.
Given a hypergraph $H = (V,E)$, with edges of cardinality at least two, and $c: V \cup E \to \Q$, the associated optimization problem is
\begin{align}
\label[problem]{prob E}
\tag{E}
\begin{split}
\max_x & \quad  \sum_{k \in V} c_k x_k + \sum_{e \in E} c_e \prod_{k \in E} x_k \\
\st & \quad x_k \in \bra{0,1} \qquad \forall k \in V.
\end{split}
\end{align}
We call this problem \emph{Explicit binary polynomial optimization,} since the polynomial objective function is given explicitly via its nonzero monomials.
This problem is strongly NP-hard \cite{GarJohSto76}, and a recent stream of research led to the discovery of several polynomially solvable classes 
\cite{dPKha17MOR,dPKha18MPA,dPKha18OS,dPKha18SIOPT,BucCraRod18,dPKhaSah20MPC,dPKha21MOR,dPKha23MPA,dPDiG21IJO,dPDiG22SODA,CapdPDiG24,dPDiG23ALG,BieMun18}.

In this paper, we propose a different way of formulating binary polynomial optimization, and show that it leads to the discovery of a vast class of tractable instances.
In this different formulation, we assume that 
the variables are partitioned into $s$ sets, and that the polynomial objective function is given in a factorized form, written as a sum of products of linear functions, each one involving only variables in one set of the partition.
Formally, this is the optimization problem
\begin{align}
\label[problem]{prob F}
\tag{F}
\begin{split}
\max_{x^1, \dots, x^s} & \quad \sum_{I \in \I} \prod_{j \in I} \inner{c^{I,j}}{x^j} \\
\st & \quad x^j \in \bra{0,1}^{n_j} \qquad \forall j \in [s],
\end{split}
\end{align}
where we denote by $\inner{a}{b}$ the inner product of two vectors $a,b$ of the same dimension.
The data of the problem consists of positive integers $s$ and $n_j$, for $j \in [s]$, 
a nonempty family $\I$ of subsets of $[s]$,
and vectors $c^{I,j} \in \Q^{n_j}$ 
for $I \in \I$, $j \in I$.
We call this problem \emph{Factorized binary polynomial optimization.}

Reading this paper, the reader should always think to the following two parameters as being fixed numbers: $\max\bra{\abs{e} : e \in E}$ in \cref{prob E}, and $s$ in \cref{prob F}.
On the one hand, this is the setting in most applications of practical interest.
On the other hand, by far most research so far is devoted to the cases in which these two parameters are equal to $2$ \cite{KocHaoGloLewLuWanWan14,Pun22}, which are already compelling and strongly NP-hard \cite{GarJohSto76,Pun22}.

While quite different on the surface, we will show that \cref{prob E,prob F} are, in fact, equivalent: we can reformulate \cref{prob E} as \cref{prob F}, and vice versa.
A key advantage of \cref{prob F} is that it can reveal some inner sparsity of the problem that is hidden in \cref{prob E}.
For example, classes of \cref{prob E} with $\abs{E}$ as large as $\prod_{j \in [s]} n_j$, can be written in the form of \cref{prob F}, with $\abs{\I} = 1$.
In this paper, we show how we can exploit this inner sparsity of the input data to design an efficient algorithm.
The main result of this paper is an algorithm that solves in strongly polynomial time \cref{prob F}, when $\abs{\I}$ is a fixed number.

\begin{theorem}
\label{th main com}
\cref{prob F} can be solved in strongly polynomial time, for any fixed $s$ and $\abs{\I}$.
\end{theorem}

It is interesting to observe that, in \cref{th main com}, it is fundamental that, in each product in the objective function of \cref{prob F}, there is at most one linear function for each set of the partition.
In fact, it is known that it is NP-hard to maximize the product of just two linear functions, plus a third linear function, over all binary points \cite{HamHanParRad02,CelKliMey06}.
As a corollary to \cref{th main com}, we obtain that also the slightly more general optimization problem obtained from \cref{prob F} by replacing each linear function in the objective with an affine function, can be solved in strongly polynomial time under the same assumptions (see \cref{cor affine}).

Our algorithm for \cref{prob F} is combinatorial in nature, and is based on the construction of hyperplane arrangements.
It is well understood that enumerating cells of hyperplane arrangements, or equivalently vertices of a zonotope, can help in the design of efficient algorithms with theoretical guarantees (see, e.g., \cite{AllFukLieSte01,FerFukLie03,dPDeyWei20SIOPT,dP23aMPA}).
Most of these algorithms construct only one arrangement in their execution.
A key feature of our algorithm is that we construct, recursively, $s-1$ arrangements on top of each other.
In the first iteration, we construct an arrangement of hyperplanes corresponding to the vector of variables $x^s$.
These hyperplanes live in a carefully chosen low-dimensional ``dual'' space, corresponding to the products in the objective function that involve variables $x^s$, and at least another $x^j$.
For each cell of the arrangement, we are able to construct a corresponding partial solution $\bar x^s \in \bra{0,1}^{n_s}$, and a child instance of \cref{prob F} only in variables $x^1,x^2,\dots,x^{s-1}$.
We then recurse on each child instance.
This gives rise to a compact tree of subproblems, and at least one of the leaves will correspond to an optimal solution of \cref{prob F}.
Our technique results in a quite elegant algorithm that is simple to state and implement, and only relies on the construction of hyperplane arrangements.

\smallskip

A consequence of the equivalence of \cref{prob E,prob F}, is that \cref{prob F} inherits the applications of \cref{prob E}.
To demonstrate the applicability of \cref{th main com}, we consider the \emph{rank-$t$ binary tensor factorization problem} and the \emph{rank-$t$ Boolean tensor factorization problem.}
In these two problems, the goal is to find $t$ rank-$1$ binary tensors whose sum, or disjunction (respectively), best approximates a given tensor of order $s$.
%
These problems arise in mining discrete patterns in data and have plenty of applications \cite{ErdMie13,ParOhKan17,RukHolYau18,WanChaZhaLiCaoZha20}.
Their solutions provide a very useful tool for analyzing tensors to discover latent factors from them, and produce more interpretable and sparser results than normal factorization methods~\cite{Mie11}.
Most existing methods to tackle these problems rely on heuristics, and hence do not provide any guarantee on the quality of the solution~\cite{Mie11,ErdMie13,BelGloVyc13,ParOhKan17,RukHolYau18}.
The reader should always think to the parameters $s$ and $t$, in these problems, as being fixed numbers:
the case $s=2$, $t=1$ is the \emph{rank-1 Boolean matrix factorization problem,} which is already NP-hard \cite{GilVav18}, well-studied, and with many applications (see, e.g., \cite{GilGli11,KoyGraRam05,KoyGraRam06,LuVaiAtlShiJia11,SheJiYe09}).
We will see that \cref{th main com} implies the following tractability results.

\begin{corollary}
\label{cor BTF rank}
The rank-$t$ binary tensor factorization problem on a tensor of order $s$ can be solved in strongly polynomial time if $s,t$ are fixed, the input tensor has fixed rank, and a rank factorization is given.
\end{corollary}

\begin{corollary}
\label{cor r1BTF rank}
The rank-1 Boolean tensor factorization problem on a tensor of order $s$ can be solved in strongly polynomial time if $s$ is fixed, the input tensor has fixed rank, and a rank factorization is given.
\end{corollary}

In fact, we do not even need a rank factorization of the input tensor, and any factorization of the tensor as the sum of a fixed number of rank-1 tensors is sufficient (see \cref{cor BTF factorization,cor r1BTF factorization}).
For the rank-1 Boolean matrix factorization problem, we obtain:

\begin{corollary}
\label{cor r1BMF rank}
The rank-1 Boolean matrix factorization problem can be solved in strongly polynomial time if the input matrix has fixed rank.
\end{corollary}


\smallskip

While the main contribution of this paper is \cref{th main com} and its consequences, 
in this work we also lay the foundations for a theoretical study of \cref{prob F}, since, to the best of our knowledge, this problem has not been considered before.
In \cref{sec hardness}, we see that \cref{prob F} is strongly NP-hard even in very restrictive settings.
In \cref{sec vs}, we discuss the relationship of \cref{prob E,prob F}:
we show that the two problems are equivalent, that \cref{prob F} can reveal some inner sparsity of the problem that is hidden in \cref{prob E}, and we provide a link with the concept of tensor factorization and rank of a tensor. 
\cref{sec BTF} is devoted to the binary tensor factorization problem and the Boolean tensor factorization problem.
Our algorithm, the proof of \cref{th main com} and its corollaries, are presented in \cref{sec proofs}.

In \cref{sec comparison}, we compare our \cref{th main com} with known tractability results.
The focus of this paper is on polynomial objective functions of degree three or more.
In this case, we are not aware of any other known polynomially solvable class of \cref{prob F}, and 
we show that our result does not follow from known tractable classes of \cref{prob E}.
Interestingly, our result also improves on the state-of-the-art for quadratic objective functions.
Consider the optimization problem
\begin{align}
\label[problem]{prob Q}
\tag{Q}
\begin{split}
\max_{x^1,\dots,x^s} & \quad \sum_{i,j \in [s], \ i<j}{x^i}^\transp Q^{i,j} x^j + \sum_{j \in [s]} {c^j}^\transp x^j \\
\st & \quad x^j \in \bra{0,1}^{n_j} \qquad \forall j \in [s],
\end{split}
\end{align}
where $Q^{i,j} \in \Q^{n_i \times n_j}$, for every $i,j \in [s]$ with $i<j$, and $c^j \in \Q^{n_j}$, for every $j \in [s]$.
As a corollary to \cref{th main com}, we obtain the following result:
\begin{corollary}
\label{cor main quad}
\cref{prob Q} can be solved in strongly polynomial time, provided that $s$ is fixed, and the rank of each matrix $Q^{i,j}$, for $i,j \in [s]$ with $i<j$, is fixed.
\end{corollary}
It is known that \cref{prob Q} can be solved in polynomial time if $s=2$ and the rank of $Q^{1,2}$ is fixed \cite{PunSriKar15a,HlaCerRad21,Pun22}.
Therefore, \cref{cor main quad} extends previously known tractable classes, both in terms of generality and running time.
Furthermore, when $s=2$, our algorithm constructs only one hyperplane arrangement, and our work results in significantly cleaner and shorter algorithm and arguments.

We remark that the main emphasis of this paper lies on the theoretical computational complexity of the algorithms presented. 
We do not refer to practically efficient implementations of the algorithms, which we believe should be studied in the future.

\section{NP-hardness}
\label{sec hardness}

In the next result, we show that \cref{prob F} is strongly NP-hard even in very restrictive settings.
Our reduction is inspired from the one in the proof of Theorem~10.2 in \cite{Pun22}, where the authors show that \cref{prob Q} is strongly NP-hard for $s=2$.
The similarity lies in the introduction of a copy of the original variables, and the addition of a penalty in the objective function to have the original variables match the copy, in an optimal solution.
The main differences are:
1) Our reduction is form \emph{simple max cut,} rather than from \emph{binary quadratic optimization,} which allows us to show the hardness of the problem in very restrictive settings;
2) Our reduction is to \cref{prob F} with $s=2$, whose objective function is given in a factorized form.

\begin{proposition}
\label{prop FMPO hard}
\cref{prob F} is strongly NP-hard even if $s=2$, $n_1 = n_2$, $\abs{\I}=n_1 + 2$, and $c^{I,j}_k$ is integer and bounded by $n_1$ in absolute value, for every $I \in \I$, $j \in I$, $k \in [n_j]$.
\end{proposition}

\begin{prf}
Our reduction is from simple max cut, which is strongly NP-hard \cite{GarJohSto76}.
It is well-known that the simple max cut problem on a graph $G = (V,E)$ with $V=[n]$ can be formulated as
\begin{align}
\label[problem]{prob max-cut}
\tag{MC}
\begin{split}
\max_{x} & \quad \sum_{\bra{i,j} \in E, \ i < j} \pare{x_i + x_j - 2 x_i x_j} \\
\st & \quad x \in \bra{0,1}^{n}.
\end{split}
\end{align}
Now consider the optimization problem
\begin{align}
\label[problem]{prob max-cut'}
\tag{MC'}
\begin{split}
\max_{x,y} & \quad \sum_{\bra{i,j} \in E, \ i < j} \pare{x_i + y_j - 2 x_i y_j} 
+ n \sum_{i \in [n]} \pare{ 2  x_i y_i - x_i - y_i } \\
\st & \quad x,y \in \bra{0,1}^{n}.
\end{split}
\end{align}
We claim that \cref{prob max-cut} can be solved by \cref{prob max-cut'}, in the sense that, given an optimal solution $(x,y)$ to \cref{prob max-cut'}, then $x$ is optimal to \cref{prob max-cut}.

Denote by $\obj(x)$ the objective value of $x$ in \cref{prob max-cut}
and by $\obj'(x,y)$ the objective value of $(x,y)$ in \cref{prob max-cut'}.
It is simple to check that, given $x,y \in \bra{0,1}^{n}$ with $x=y$, then $\obj'(x,y) = \obj(x)$.
It then suffices to show that if $(x,y)$ is optimal to \cref{prob max-cut'}, then $x=y$.
Let $(x,y)$ be optimal to \cref{prob max-cut'} and assume, for a contradiction, that there is at least one index $\ell \in [n]$ such that $x_\ell \neq y_\ell$.
Let $\tilde y$ be obtained from $y$ by flipping component $\ell$, so that $\tilde y_\ell = x_\ell$. 
We show that the objective value of $(x,\tilde y)$ is strictly larger than the objective value of $(x,y)$, which contradicts the optimality of $(x,y)$:
\begin{align*}
\obj'(x, y) - \obj'(x, \tilde y) 
& = \sum_{\bra{i,j} \in E, \ i < j} \pare{1 - 2 x_i} \pare{y_j - \tilde y_j}
+ n \sum_{i \in [n]} \pare{2 x_i - 1} \pare{y_i - \tilde y_i} \\
& = \sum_{i : \bra{i,\ell} \in E, \ i < \ell} \pare{1 - 2 x_i} \pare{y_\ell - \tilde y_\ell}
+ n \pare{2 x_\ell - 1} \pare{y_\ell - \tilde y_\ell} \\
& = \sum_{i : \bra{i,\ell} \in E, \ i < \ell} \pare{1 - 2 x_i} \pare{1- 2 x_\ell}
+ n \pare{2 x_\ell - 1} \pare{1- 2 x_\ell} \\
& = \absL{\bra{i \in [\ell-1] : \bra{i,\ell} \in E}} - n \\
& \le (n-1) - n = -1.
\end{align*}

To complete the proof, we write \cref{prob max-cut'} in the form of \cref{prob F} with $s=2$, $n_1=n_2=n$, $\abs{\I}=n + 2$, $x^1 = x$, and $x^2 = y$.
We rewrite the objective function of \cref{prob max-cut'} as the sum of the following three functions:
\begin{align*}
f(x) & := \sum_{\bra{i,j} \in E, \ i < j} x_i
- n \sum_{i \in [n]} x_i \\
& \phantom{:}= \sum_{i \in [n]} x_i \pare{ \abs{\bra{j : \bra{i,j} \in E, \ i < j}} - n}, \\
g(y) & := \sum_{\bra{i,j} \in E, \ i < j} y_j
- n \sum_{j \in [n]} y_j \\
& \phantom{:}= \sum_{j \in [n]} y_j \pare{ \abs{\bra{i : \bra{i,j} \in E, \ i < j}} - n}, \\
h(x,y) & := -2 \sum_{\bra{i,j} \in E, \ i < j} x_i y_j
+ 2 n \sum_{i \in [n]} x_i y_i \\
& \phantom{:}= \sum_{i \in [n]} \underbrace{2x_i}_{=:h_i(x)} \underbrace{\pare{n y_i - \sum_{j : \bra{i,j} \in E, \ i < j} y_j}}_{=:h_i(y)}.
\end{align*}
Now, the function $f(x)$ corresponds to a set $I=\bra{1}$ with $c^{I,1}$ the vector of coefficients of the linear function $f(x)$.
The function $g(y)$ corresponds to a set $I=\bra{2}$ with $c^{I,2}$ the vector of coefficients of the linear function $g(y)$.
For every $i \in [n]$, the product $h_i(x)h_i(y)$ in $h(x,y)$ corresponds to a set $I=\bra{1,2}$ with $c^{I,1}$ and $c^{I,2}$ the vectors of coefficients of the linear functions $h_i(x)$ and $h_i(y)$, respectively.
It is simple to check that all these vectors have integer components bounded by $n$ in absolute value.
\end{prf}

\section{Binary polynomial optimization: explicit vs factorized}
\label{sec vs}

In this section we compare \cref{prob E} and \cref{prob F}.

\subsection{Equivalency of \cref{prob E,prob F}}
\label{sec equivalency E F}

Our first goal is to establish the equivalency of \cref{prob E,prob F}.
While it is clear that we can directly reformulate \cref{prob F} as \cref{prob E}, by expanding all products in the objective function, the opposite direction is less obvious.
To reformulate \cref{prob E} as \cref{prob F}, we expand the technique that we used in the proof of \cref{prop FMPO hard}, from quadratics to higher degree polynomials: the constructed instance of \cref{prob F} features $s$ copies of the original variables, where $s$ denotes the degree of the objective function of \cref{prob E}, and a penalty in the objective function forces all the copies to match, in an optimal solution.

In the remainder of the paper, we say that a hypergraph $H$ is \emph{$s$-partite,} if its node set can be partitioned into $s$ sets, called the \emph{sides} of $H$, such that every edge contains at most one node from each side.

\begin{proposition}
\label{prop reduction E F}
\cref{prob E} on a hypergraph $H=(V,E)$ can be reformulated, in strongly polynomial time, as \cref{prob F} with $s := \max\bra{\abs{e} : e \in E}$ and $\abs{\I} \le \abs{V} + \abs{E} + s$.
\end{proposition}

\begin{prf}
Consider an instance of \cref{prob E}, and let $\obj(x)$ be its objective function.
Construct a new instance of \cref{prob E} as follows:
1) The instance has variables $x^1,x^2, \dots,x^s$, where $s := \max\bra{\abs{e} : e \in E}$, and each $x^j$, for $j \in [s]$, is a copy of the original vector of variables $x$ of the original instance.
2) Let $f(x^1,x^2, \dots,x^s)$ be obtained from $\obj(x)$ by replacing, for every $e \in E$, each $x_k$ with one of its copies $x^j_k$, for $j \in [s]$, so that the product $\prod_{k \in E} x_k$ contains at most one variable from each $x^j$, $j \in [s]$.
3) The objective function of the new instance is then defined by
\begin{align*}
\obj'(x^1,x^2, \dots,x^s) := f(x^1,x^2, \dots,x^s) + M \sum_{k \in V} \pare{s \prod_{j \in [s]} x^j_k - \sum_{j \in [s]} x^j_k},
\end{align*}
where $M := \sum_{e \in E} \abs{c_e}+1$.

Note that the hypergraph $H'=(V',E')$ associated with the new instance of \cref{prob E} is $s$-partite, and is constructed as follows: 
1) $V'$ contains $s$ copies of $V$, which we denote by $V^j$, for $j \in [s]$.
2) For every edge $e \in E$, we have an edge $e' \in E'$ obtained from $e$ by replacing each node with one if its copies, so that $e'$ contains at most one node from each $V^j$, for $j \in [s]$;
Furthermore, for every $k \in V$, we have an edge in $E'$ that contains all $s$ copies of $k$.
Therefore, $\abs{E'} \le \abs{V} + \abs{E}$.

Let $(x^1,x^2, \dots,x^s)$ be an optimal solution to the new instance.
We show that $x^1$ is optimal to the original instance.
Define, for every $k \in V$, the penalty function 
\begin{align*}
p_k (x^1_k,x^2_k,\dots,x^s_k) := 
s \prod_{j \in [s]} x^j_k 
- \sum_{j \in [s]} x^j_k.
\end{align*}
Given $(x^1,x^2,\dots,x^s)$ binary with $x^1 = x^2 = \dots = x^s$, we have $p_k (x^1_k,x^2_k,\dots,x^s_k) = 0$ for every $k \in V$, thus $\obj'(x^1,x^2,\dots,x^s) = \obj(x^1)$.
It then suffices to show that if $(x^1,x^2,\dots,x^s)$ is optimal to the new instance, then $x^1 = x^2 = \dots = x^s$.
Let $(x^1,x^2,\dots,x^s)$ be optimal to the new instance and assume, for a contradiction, that there is at least one $\ell \in V$ such that $x^1_\ell = x^2_\ell = \dots = x^s_\ell$ does not hold.
Let $(\tilde x^1,\tilde x^2,\dots,\tilde x^s)$ be obtained from $(x^1,x^2,\dots,x^s)$ by setting all components $\tilde x^1_\ell = \tilde x^2_\ell = \dots = \tilde x^s_\ell$ to the same binary value.
We show that $\obj'(\tilde x^1,\tilde x^2,\dots,\tilde x^s)$ is strictly larger than $\obj'(x^1,x^2,\dots,x^s)$, which contradicts the optimality of $(x^1,x^2,\dots,x^s)$.
First, by definition of $M$, 
\begin{align*}
f(x^1,x^2,\dots,x^s) - f(\tilde x^1,\tilde x^2,\dots,\tilde x^s) < M.
\end{align*}
Next, 
we have 
\begin{align*}
& p_\ell(x^1_\ell,x^2_\ell,\dots,x^s_\ell) - p_\ell(\tilde x^1_\ell,\tilde x^2_\ell,\dots,\tilde x^s_\ell) \le -1 \\
& p_k(x^1_k,x^2_k,\dots,x^s_k) - p_k(\tilde x^1_k,\tilde x^2_k,\dots,\tilde x^s_k) = 0 \qquad \forall k \in V \setminus \bra{\ell}.
\end{align*}
Hence,
\begin{align*}
& \obj'(x^1,x^2,\dots,x^s) - \obj'(\tilde x^1,\tilde x^2,\dots,\tilde x^s) 
 < M - M = 0.
\end{align*}

To complete the proof, we write the new instance in the form of \cref{prob F} with $s$ as above, $n_j = \abs{V}$ for every $j \in [s]$, and $\abs{\I} \le \abs{V} + \abs{E} + s$.
For every $j \in [s]$, the part of the objective function of the new instance that is linear in variables $x^j$ corresponds to a set $I=\bra{j}$ with $c^{I,j}$ the vector of coefficients of the variables $x^j$.
Next, consider a single nonlinear monomial in the objective function of the new instance, which we can write as $c \prod_{j \in I} x^j_{k_j}$, for some $I \subseteq [s]$ and $k_j \in [n_j]$, for $j \in [s]$.
This monomial can be written as $\prod_{j \in I} \inner{c^{I,j}}{x^j}$, where
the $c^{I,j}_{k_j}$, for every $j \in I$, are set in any way so that their product is $c$, and 
every other component of $c^{I,j}$, for every $j \in I$, is set to zero.
\end{prf}

\subsection{Tensor factorizations and rank}

In \cref{sec equivalency E F}, we saw that \cref{prob E,prob F} can be polynomially reduced into each other.
However, in some cases, one formulation can be significantly more succinct than the other.
In the next example we see that, starting from an instance of \cref{prob F} with $\abs{\I}=1$, and expanding all products in the objective function, it is possible to obtain an instance of \cref{prob E} with $\prod_{j \in [s]} n_j$ monomials with nonzero coefficients.
In particular, the hypergraph associated with the obtained instance of \cref{prob E} is $s$-partite, $s$-uniform, and contains all possible $\prod_{j \in [s]} n_j$ edges.
Recall that a hypergraph is \emph{$s$-uniform} if all its edges have cardinality $s$.

\begin{example}
\label{ex cu}
Consider the special case of \cref{prob F} with $\I = \bra{I}$, $I = [s]$, and with each vector $c^{I,j} \in \Q^{n_j}$, for every $j \in [s]$, having all nonzero components.
We expand all products in the objective function and obtain:
\begin{align*}
\quad \prod_{j \in [s]} \pare{\sum_{k_j \in [n_j]}{c^{I,j}_{k_j} x^j_{k_j}}} =
\sum_{k_1 \in [n_1]} \cdots \sum_{k_s \in [n_s]}
c^{I,1}_{k_1} c^{I,2}_{k_2} \cdots c^{I,s}_{k_s} \cdot x^1_{k_1} x^2_{k_2} \cdots x^s_{k_s}.
\end{align*}
The hypergraph associated with the obtained instance of \cref{prob E} is $s$-partite and $s$-uniform.
Since each product $c^{I,1}_{k_1} c^{I,2}_{k_2} \cdots c^{I,s}_{k_s}$, for $k_1 \in [n_1], k_2 \in [n_2], \dots, k_s \in [n_s]$, is nonzero, and no cancellation is possible, it contains all possible $\prod_{j \in [s]} n_j$ edges.
\end{example}

\cref{ex cu} highlights how \cref{prob F} can provide a much sparser formulation with respect to \cref{prob E}, revealing some inner sparsity of the problem that is hidden in \cref{prob E}.
Our next goal is to connect this phenomenon to the concept of tensor factorizations and rank of a tensor.

To avoid unnecessary complicated notation, we mainly focus on the special cases of \cref{prob E,prob F} that we implicitly considered in \cref{ex cu}.
For \cref{prob E}, we focus on the special case in which the associated hypergraph $H$ is $s$-partite, with sides $V^j$, for $j \in [s]$, and $s$-uniform.
As a result, $H$ has at most $\prod_{j \in [s]} n_j$ edges, where $n_j := \abs{V^j}$, for every $j \in [s]$.
We can write each edge $e \in E$ in the form $e=\pare{k_1,k_2,\dots,k_s}$, for some $k_1 \in [n_1]$, $k_2 \in [n_2]$, \dots, $k_s \in [n_s]$.
We can then write this special case of \cref{prob E} in the form
\begin{align}
\label[problem]{prob E-U}
\tag{E-U}
\begin{split}
\max_{x^1, \dots, x^s} & \quad \sum_{\pare{k_1,\dots,k_s} \in E} c_{k_1,\dots,k_s} \prod_{j \in [s]} x^j_{k_j} \\
\st & \quad x^j \in \bra{0,1}^{n_j} \qquad \forall j \in [s].
\end{split}
\end{align}
Note that the degree of each monomial in the objective function is $s$.
Furthermore, we can define the associated $n_1 \times n_2 \times \cdots \times n_s$ tensor of the coefficients of the objective function $C := \pare{c_{k_1,\dots,k_s}}$, with the understanding that $c_{k_1,\dots,k_s}$ is set to zero, if $\pare{k_1,\dots,k_s} \notin E$.
%
On the other hand, in the special case of \cref{prob F} that we consider, we have that all sets $I \in \I$ coincide with the set $[s]$.
We can then set $r := \abs{\I}$ in \cref{prob F} and write this special case of \cref{prob F} in the form:
\begin{align}
\label[problem]{prob F-U}
\tag{F-U}
\begin{split}
\max_{x^1, \dots, x^s} & \quad \sum_{i \in [r]} \prod_{j \in [s]} \inner{c^{i,j}}{x^j} \\
\st & \quad x^j \in \bra{0,1}^{n_j} \qquad \forall j \in [s].
\end{split}
\end{align}

Next, we show that \cref{prob F-U} provides a much sparser formulation with respect to \cref{prob E-U}, when the tensor $C$ of the coefficients of the objective function of \cref{prob E-U} can be written as a sum of few rank-1 tensors:
\begin{align}
\label{eq factorization}
\sum_{i \in [r]} (c^{i,1} \otimes c^{i,2} \otimes \cdots \otimes c^{i,s}),
\end{align}
where $\otimes$ denotes the vector outer product, and
where $c^{i,j} \in \Q^{n_j}$, for every $i \in [r]$ and $j \in[s]$.
In this case, we say that \eqref{eq factorization} is a \emph{factorization} of $C$.

\begin{observation}
\label{obs rank to factorization}
Assume that the tensor $C$ of the coefficients of the objective function of \cref{prob E-U} has a factorization of the form \eqref{eq factorization}.
Then, \cref{prob E-U} can be reformulated as \cref{prob F-U}, where $r$ and the $c^{i,j}$ are the ones from \eqref{eq factorization}.
\end{observation}

\begin{prf}
Let $C = \pare{c_{k_1,\dots,k_s}}$ be the $n_1 \times n_2 \times \cdots \times n_s$ tensor of the coefficients of the objective function of \cref{prob E-U}.
Since it has a factorization of the form \eqref{eq factorization}, 
the entries of $C$ are
\begin{align*} 
c_{k_1,k_2,\dots,k_s} = \sum_{i \in [r]} 
\prod_{j \in [s]} c^{i,j}_{k_j}
\qquad \forall k_1 \in [n_1], \ k_2 \in [n_2], \dots, k_s \in [n_s].
\end{align*} 
The objective function of \cref{prob E-U} can then be written in the form of the objective function  \cref{prob F-U} as follows:
\begin{align*} 
\sum_{\pare{k_1,\dots,k_s} \in E} c_{k_1,\dots,k_s} \prod_{j \in [s]} x^j_{k_j}
& = \sum_{k_1 \in [n_1]} \cdots \sum_{k_s \in [n_s]} c_{k_1,\dots,k_s} \prod_{j \in [s]} x^j_{k_j} \\
& = \sum_{k_1 \in [n_1]} \cdots \sum_{k_s \in [n_s]} \pare{\sum_{i \in [r]} \prod_{j \in [s]} c^{i,j}_{k_j}} \prod_{j \in [s]} x^j_{k_j} \\
& = \sum_{i \in [r]} \sum_{k_1 \in [n_1]} \cdots \sum_{k_s \in [n_s]} \prod_{j \in [s]} c^{i,j}_{k_j} x^j_{k_j} \\
& = \sum_{i \in [r]} \prod_{j \in [s]} \pare{\sum_{k_j \in [n_j]} c^{i,j}_{k_j} x^j_{k_j}} \\
& = \sum_{i \in [r]} \prod_{j \in [s]} \inner{c^{i,j}}{x^j}.
\end{align*} 
\end{prf}

\cref{obs rank to factorization} allows us to draw a connection with the concept of tensor rank.
The \emph{rank} of an $n_1 \times n_2 \times \cdots \times n_s$ tensor $C$ is the minimum number $r$ such that $C$ has a factorization of the form \eqref{eq factorization}.
In this case, we say that \eqref{eq factorization} is a \emph{rank factorization} of $C$.
It is well-known that the rank of an $n_1 \times n_2 \times \cdots \times n_s$ tensor $C$ with $s \ge 2$ is less than or equal to the product of the $s-1$ smallest numbers among $n_1, n_2, \dots,n_s$.
In the case $s=2$, the definition of rank and of rank factorization of a tensor reduce to the well-known definitions for a matrix.
In fact, we have
\begin{align*}
\sum_{i \in [r]} (c^{i,1} \otimes c^{i,2}) = AB^\transp,
\end{align*}
where $A$ is the $n_1 \times r$ matrix with column $i$ equal to $c^{i,1}$,
and $B$ is the $n_2 \times r$ matrix with column $i$ equal to $c^{i,2}$, for $i \in [r]$.
While it is possible to determine the rank of a matrix in polynomial time, for example via Gaussian elimination, determining the rank of a tensor is NP-hard, even for a tensor of order 3.
This follows from \cite{Has90}, which details how to encode any given 3SAT Boolean formula in $n$ variables and $m$ clauses as an $(n+2m+2) \times 3n \times (3n+m)$ tensor $C$ with the property that the 3SAT formula is satisfiable if and only if the rank of $C$ is at most $4n + 2m$.

The next result characterizes the best possible sparsity of \cref{prob F-U} in terms of the rank of the tensor of the coefficients of the objective function of \cref{prob E-U}.

\begin{observation}
\label{obs tight rank}
\cref{prob E-U} can be reformulated as \cref{prob F-U}, where the number $r$ in \cref{prob F-U} equals the rank of the tensor $C$ of the coefficients of the objective function of \cref{prob E-U}.
This is best possible, in the sense that it is not possible to write the objective function of \cref{prob E-U} in the form of the objective function of \cref{prob F-U} with a strictly smaller $r$.
\end{observation}

\begin{prf}
The first part of the statement follows directly from the definition of rank and \cref{obs rank to factorization}.
In the remainder of the proof we prove the second part of the statement.
Assume that the objective function of \cref{prob E-U} can be written in the form of the objective function of \cref{prob F-U}:
\begin{align*} 
\sum_{i \in [r]} \prod_{j \in [s]} \inner{c^{i,j}}{x^j},
\end{align*} 
where $c^{i,j} \in \Q^{n_j}$, for every $i \in [r]$ and $j \in[s]$.
Let $C = \pare{c_{k_1,\dots,k_s}}$ be the $n_1 \times n_2 \times \cdots \times n_s$ tensor of the coefficients of the objective function of \cref{prob E-U}.
Following the same derivation in the proof of \cref{obs rank to factorization}, but in the reverse order, we obtain that the entries of $C$ are
\begin{align*} 
c_{k_1,k_2,\dots,k_s} = \sum_{i \in [r]} 
\prod_{j \in [s]} c^{i,j}_{k_j}
\qquad \forall k_1 \in [n_1], \ k_2 \in [n_2], \dots, k_s \in [n_s].
\end{align*} 
Hence, $C$ has a factorization of the form \eqref{eq factorization}.
From the definition of rank, we have that $r$ is greater than or equal to the rank of $C$.
\end{prf}

The discussion above about the complexity of determining the rank of a tensor, highlights a key difference in the complexity of the reductions in \cref{prop reduction E F,obs rank to factorization,obs tight rank}: 
The reformulation in \cref{prop reduction E F} can be constructed in strongly polynomial time.
On the other hand, the reformulation in \cref{obs rank to factorization} can be obtained in strongly polynomial time if the factorization \eqref{eq factorization} of $C$ is given.
Furthermore, the reformulation in \cref{obs tight rank} can be obtained in strongly polynomial time if a rank factorization of $C$ is given, and it cannot be obtained in polynomial time, in general, unless P=NP.

\smallskip

\cref{obs rank to factorization,obs tight rank} are stated in terms of \cref{prob E-U,prob F-U}.
Similar results can be also obtained if we consider, instead of \cref{prob E-U}, \cref{prob E} in which the associated hypergraph is $s$-partite with sides $V^j$, for $j \in I$.
In other words, \cref{prob E-U} without the assumption that the hypergraph is $s$-uniform.
A key difference, which results in a significantly more complicated notation, is that, in this more general problem, we can no longer encode the coefficients of the objective function with only one tensor.
Instead, we need at most $2^s - 1$ tensors: one for each subset $I$ of $[s]$ such that there is an edge that contains one node in each $V^j$, for $j \in I$.
In these results, we then need to replace \cref{prob F-U} with \cref{prob F}.
For this pair of more general problems, one can obtain the following results, which correspond to \cref{obs rank to factorization,obs tight rank}.

\begin{observation}
\label{obs rank to factorization full}
Consider \cref{prob E} on an $s$-partite hypergraph.
Assume that all the tensors of the coefficients of the objective function have a factorization of the form \eqref{eq factorization}.
Then, this problem can be reformulated as \cref{prob F}, where $\abs{\I}$ equals the sum of all the $r$ in \eqref{eq factorization}, and the $c^{I,j}$ are the ones from \eqref{eq factorization}.
\end{observation}

\begin{observation}
\label{obs tight rank full}
Consider \cref{prob E} on an $s$-partite hypergraph.
This problem can be reformulated as \cref{prob F}, where $\abs{\I}$ in \cref{prob F} equals the sum of the ranks of the tensors of the coefficients of the objective function of \cref{prob E}.
This is best possible, in the sense that it is not possible to write the objective function of \cref{prob E} in the form of the objective function of \cref{prob F} with a strictly smaller $\abs{\I}$.
\end{observation}

\cref{obs rank to factorization full,obs tight rank full} follow by applying \cref{obs rank to factorization,obs tight rank}, respectively, to each tensor of the coefficients of the objective function of the considered \cref{prob E}.



\section{Binary tensor factorization and Boolean tensor factorization}
\label{sec BTF}


To formally introduce the binary tensor factorization problem, and the Boolean tensor factorization problem, we first define the binary rank and the Boolean rank of a tensor.
The \emph{binary rank} of an $n_1 \times n_2 \times \cdots \times n_s$ tensor $C$ is the minimum number $r$ such that $C$ has a factorization of the form \eqref{eq factorization}, where $c^{i,j} \in \bra{0,1}^{n_j}$, for every $i \in [r]$ and $j \in[s]$.
On the other hand, the \emph{Boolean rank} of an $n_1 \times n_2 \times \cdots \times n_s$ tensor $C$ is the minimum number $r$ such that $C$ has a factorization of the form
\begin{align}
\label{eq Boolean factorization}
\bigvee_{i \in [r]} (c^{i,1} \otimes c^{i,2} \otimes \cdots \otimes c^{i,s}),
\end{align}
where $\vee$ denotes the component-wise ``or'' operation, and 
where $c^{i,j} \in \bra{0,1}^{n_j}$, for every $i \in [r]$ and $j \in[s]$.
We refer the interested reader to \cite{HavPar22} for more informations on the different notions of rank.
%
%
In the \emph{rank-$t$ binary tensor factorization problem,} we are given an $n_1 \times n_2 \times \cdots \times n_s$ tensor $A=\pare{a_{k_1,\dots,k_s}}$ and an integer $t$.
The goal is to find an $n_1 \times n_2 \times \cdots \times n_s$ tensor $B=\pare{b_{k_1,\dots,k_s}}$ of binary rank $t$
that minimizes the objective function
\begin{align}
\label{eq tensor obj}
\sum_{k_1 \in [n_1]} \cdots \sum_{k_s \in [n_s]} \pare{a_{k_1,\dots,k_s}-b_{k_1,\dots,k_s}}^2.
\end{align}
Also in the \emph{rank-$t$ Boolean tensor factorization problem,} we are given an $n_1 \times n_2 \times \cdots \times n_s$ tensor $A=\pare{a_{k_1,\dots,k_s}}$ and an integer $t$.
However, the goal is to find an $n_1 \times n_2 \times \cdots \times n_s$ tensor $B=\pare{b_{k_1,\dots,k_s}}$ of Boolean rank $t$ that minimizes the objective function \eqref{eq tensor obj}.
%
In the special case $t=1$ the two problems coincide.
The case $s=2$, $t=1$ is the \emph{rank-1 Boolean matrix factorization problem.} 
It is shown in \cite{Pun22} that the latter problem can be formulated as \cref{prob Q} with $s=2$.
Next, we extend this result to tensors of any order, and to the more general rank-$t$ binary tensor factorization problem.


\begin{proposition}
\label{prop BTF E}
The rank-$t$ binary tensor factorization problem on an $n_1 \times n_2 \times \cdots \times n_s$ tensor $A=\pare{a_{k_1,\dots,k_s}}$ can be formulated as the following \cref{prob E} on an $st$-partite hypergraph:
\begin{align}
\label[problem]{prob BTF}
\tag{BTF}
\begin{split}
\max_{x^{i,j}} & \quad \sum_{k_1 \in [n_1]} \cdots \sum_{k_s \in [n_s]}
\sbra{\pare{2 a_{k_1,\dots,k_s} - 1} \sum_{i \in [t]} \prod_{j \in [s]} x^{i,j}_{k_j}  
- 2 \sum_{i,i' \in [t], \ i < i'} \prod_{j \in [s]} x^{i,j}_{k_j} x^{i',j}_{k_j}} \\
\st & \quad x^{i,j} \in \bra{0,1}^{n_j} \qquad \forall i \in [t], \ j \in [s].
\end{split}
\end{align}
\end{proposition}




\begin{prf}
Consider the binary tensor factorization problem.
For every $k_1 \in [n_1], k_2 \in [n_2], \dots, k_s \in [n_s]$ we have 
\begin{align*}
b_{k_1,\dots,k_s} 
= \sum_{i \in [t]} x^{i,1}_{k_1} x^{i,2}_{k_2} \cdots x^{i,s}_{k_s}.
\end{align*}
Since $x^{i,j}$ is binary, for every $i \in [t]$, $j \in [s]$, we have, for every $k_1 \in [n_1], k_2 \in [n_2], \dots, k_s \in [n_s]$:
\begin{align*}
& \pare{a_{k_1,\dots,k_s} -b_{k_1,\dots,k_s}}^2
= \pare{a_{k_1,\dots,k_s}-\sum_{i \in [t]} \prod_{j \in [s]} x^{i,j}_{k_j} }^2 \\
& = a_{k_1,\dots,k_s}^2 - 2 a_{k_1,\dots,k_s} \sum_{i \in [t]} \prod_{j \in [s]} x^{i,j}_{k_j} 
+ \pare{\sum_{i \in [t]} \prod_{j \in [s]} x^{i,j}_{k_j}}^2 \\
& = a_{k_1,\dots,k_s}^2 - 2 a_{k_1,\dots,k_s} \sum_{i \in [t]} \prod_{j \in [s]} x^{i,j}_{k_j}  
+ \sum_{i \in [t]} \prod_{j \in [s]} \pare{x^{i,j}_{k_j}}^2
+ 2 \sum_{i,i' \in [t], \ i < i'} \prod_{j \in [s]} x^{i,j}_{k_j} x^{i',j}_{k_j} \\
& = a_{k_1,\dots,k_s}^2 - 2 a_{k_1,\dots,k_s} \sum_{i \in [t]} \prod_{j \in [s]} x^{i,j}_{k_j}  
+ \sum_{i \in [t]} \prod_{j \in [s]} x^{i,j}_{k_j}
+ 2 \sum_{i,i' \in [t], \ i < i'} \prod_{j \in [s]} x^{i,j}_{k_j} x^{i',j}_{k_j} \\
& = a_{k_1,\dots,k_s}^2 + \pare{1 - 2 a_{k_1,\dots,k_s}} \sum_{i \in [t]} \prod_{j \in [s]} x^{i,j}_{k_j}  
+ 2 \sum_{i,i' \in [t], \ i < i'} \prod_{j \in [s]} x^{i,j}_{k_j} x^{i',j}_{k_j}.
\end{align*}
Our binary tensor factorization problem can then be formulated as \cref{prob E}.
Note that, in each monomial in the objective function, no two variables from the same vector $x^{i,j}$, for $i \in [t]$, $j \in [s]$, are ever multiplied together in the same monomial.
Therefore, this optimization problem is \cref{prob E} on an $st$-partite hypergraph, with $t(n_1+n_2+ \cdots +n_s)$ variables, and with an objective function of degree $2s$.
\end{prf}

\begin{corollary}
\label{cor r1BTF E-U}
The rank-$1$ Boolean tensor factorization problem on an $n_1 \times n_2 \times \cdots \times n_s$ tensor $A=\pare{a_{k_1,\dots,k_s}}$ can be formulated as the following \cref{prob E-U} on an $s$-partite hypergraph:
\begin{align}
\label[problem]{prob 1BTF}
\tag{BTF$_1$}
\begin{split}
\max_{x^1,\dots,x^s} & \quad \sum_{k_1 \in [n_1]} \cdots \sum_{k_s \in [n_s]}
\pare{2 a_{k_1,\dots,k_s} - 1} \prod_{j \in [s]} x^j_{k_j} \\
\st & \quad x^{j} \in \bra{0,1}^{n_j} \qquad \forall j \in [s].
\end{split}
\end{align}
\end{corollary}


\begin{prf}
In the special case $t=1$, the rank-$1$ Boolean tensor factorization problem coincides with the rank-$1$ binary tensor factorization problem.
It then follows from \cref{prop BTF E} that the rank-$1$ Boolean tensor factorization problem on an $n_1 \times n_2 \times \cdots \times n_s$ tensor $A=\pare{a_{k_1,\dots,k_s}}$ can be formulated as 
\cref{prob 1BTF},
which is of the form of \cref{prob E} on an $s$-partite hypergraph.
Note that, in each monomial in the objective function, precisely one variable from each vector $x^{j}$, for $j \in [s]$, is present.
Therefore, this optimization problem is of the form of \cref{prob E-U}.
\end{prf}

Next, we see how the binary tensor factorization problem and the rank-1 Boolean tensor factorization problem can be formulated succinctly as \cref{prob F} and \cref{prob F-U}, respectively, if a factorization of the input tensor $A$ is given.
This result will allow us to show our tractability results for the binary tensor factorization problem and its special cases, \cref{cor BTF rank,cor r1BTF rank,cor r1BMF rank}, as consequences of our main result.

\begin{proposition}
\label{prop BTF F}
Consider the rank-$t$ binary tensor factorization problem on an $n_1 \times n_2 \times \cdots \times n_s$ tensor
\begin{align*}
A = \sum_{p \in [q]} \pare{a^{p,1} \otimes a^{p,2} \otimes \cdots \otimes a^{p,s}},
\end{align*}
where $a^{p,j} \in \Q^{n_j}$, for every $p \in [q]$, $j \in [s]$.
Then, the problem can be reformulated, in strongly polynomial time, as \cref{prob F} with $s:=st$ and with $\abs{\I}:=tq + (t^2+t)/2$.
\end{proposition}

\begin{prf}
Let $A=\pare{a_{k_1,\dots,k_s}}$.
From \cref{prop BTF E}, the rank-$t$ binary tensor factorization problem can be formulated as \cref{prob BTF}.
In this problem, the coefficients of the objective function can be encoded with one tensor for each sum inside the square brackets in the objective function.

First, we consider the first sum, over all $i \in [t]$, and we consider now one specific such $i$.
The corresponding tensor of coefficients is an $n_1 \times n_2 \times \cdots \times n_s$ tensor, and the entry in position $(k_1,k_2,\dots,k_s)$ is $2 a_{k_1,\dots,k_s} - 1$.
Clearly, we have 
\begin{align*}
2 A = \sum_{p \in [q]} \pare{a^{p,1} \otimes a^{p,2} \otimes \cdots \otimes 2 a^{p,s}}.
\end{align*}
For every $j \in [s]$, let $e^j$ be the vector in $\Q^{n_j}$ with all components equal to one.
Then, the $n_1 \times n_2 \times \cdots \times n_s$ tensor with all entries equal to one can be written as
\begin{align*}
e^1 \otimes e^2 \otimes \cdots \otimes (-e^s).
\end{align*}
We can now write the above tensor of coefficients in the form
\begin{align*}
\sum_{p \in [q]} \pare{a^{p,1} \otimes a^{p,2} \otimes \cdots \otimes 2 a^{p,s}} + \pare{e^1 \otimes e^2 \otimes \cdots \otimes (-e^s)}.
\end{align*}

Next, we consider the second sum inside the square brackets in the objective function.
This sum is over all $i,i' \in [t]$ with $i < i'$, and we consider now one specific such pair $i,i'$.
The corresponding tensor of coefficients is an $n_1 \times n_1 \times n_2 \times n_2 \times \cdots \times n_s \times n_s$ tensor, and all the entries are equal to $-2$.
We can then write the above tensor of coefficients as
\begin{align*}
e^1 \otimes e^2 \otimes \cdots \otimes (-2e^s).
\end{align*}

\cref{prob BTF} is a special case of \cref{prob E}.
It then follows from \cref{obs rank to factorization full} that the rank-$t$ binary tensor factorization problem can be formulated as \cref{prob F}, with 
$s:=st$ and 
$\abs{\I} = t(q+1) + (t^2-t) / 2 = tq + (t^2+t)/2.$
\end{prf}


\begin{corollary}
\label{cor r1BTF F-U}
Consider the rank-$1$ Boolean tensor factorization problem on an $n_1 \times n_2 \times \cdots \times n_s$ tensor
\begin{align*}
A = \sum_{p \in [q]} \pare{a^{p,1} \otimes a^{p,2} \otimes \cdots \otimes a^{p,s}},
\end{align*}
where $a^{p,j} \in \Q^{n_j}$, for every $p \in [q]$, $j \in [s]$.
Then, the problem can be reformulated, in strongly polynomial time, as \cref{prob F-U} with the same $s$, $n_j$, for $j \in [s]$, and with $r:=q+1$.
\end{corollary}


\paragraph{More applications.}

We conclude this section by mentioning some additional applications for which \cref{prob E}, with an associated $s$-partite hypergraph, and \cref{prob F} provide a better fit than \cref{prob E}.
In \cite{Pun22} the author provides a number of applications of \cref{prob Q} with $s=2$.
These are the \emph{maximum weight biclique problem}, the \emph{maximal sum submatrix problem}, and the problem of finding the \emph{cut-norm of a matrix.}
\cref{prob E}, with an associated $s$-partite hypergraph, allows us to directly formulate problems that contain the ones above as special cases.
In the maximum weight biclique problem, instead of a bipartite graph, we can consider a multipartite graph with more than two sides, or even a multipartite hypergraph.
In the maximal sum submatrix problem, and in the problem of finding the cut-norm of a matrix, instead of a matrix in input, we can have a tensor of any order.
If this input tensor is given through a factorization of the form \eqref{eq factorization}, then these problems can be naturally formulated in the form of \cref{prob F-U}, due to \cref{obs rank to factorization}.
We leave the details of these formulations to the reader.

\section{Proofs of main results}
\label{sec proofs}

In this section, we prove our main results.
We begin by stating an extended version of \cref{th main com}, presented in \cref{sec intro}.
In \cref{prob F}, we denote by $m_j$, for $j \in [s]$, the number of sets $I \in \I$ that contain index $j$ and at least one index strictly smaller than $j$.

\begin{theorem}
\label{th main ext}
\cref{prob F} can be solved with 
\begin{align*}
s \theta^{s-1} n_2^{m_2} n_3^{m_3} \cdots n_s^{m_s} \poly(n_1,n_2,\dots,n_s,\abs{\I})
\end{align*}
arithmetic operations, where $\theta$ is a constant.
Furthermore, the size of the numbers produced in the course of the execution of the algorithm is polynomial in the size of the input.
In particular, \cref{prob F} can be solved in strongly polynomial time, provided that $s,m_2,\dots, m_s$ are fixed.
\end{theorem}
Clearly, \cref{th main ext} implies \cref{th main com}, since $\abs{\I}$ in \cref{prob F} is an upper bound on $m_j$, for every $j \in [s]$.
It is interesting to note that the parameter $m_1$, unlike all other $m_j$, for $j \in [s] \setminus \bra{1}$, never appears as an exponent in the running time in \cref{th main ext}, so it does not need to be fixed to obtain a polynomial time algorithm.
Furthermore, due to the symmetry of the problem, we can assume without loss of generality that $m_1$ is the largest of the $m_j$, for $j \in [s]$.
This detail is lost in \cref{th main com}.
The proof of \cref{th main ext} is presented in the next section.

\subsection{Proof of \cref{th main ext}}
\label{sec proof th}

Our algorithm is based on the construction of hyperplane arrangements.
We now introduce the proper terminology.
A finite family of hyperplanes in $\R^d$ defines a dissection of $\R^d$ into connected sets of various dimensions.
We call this dissection the \emph{arrangement} of these hyperplanes.
The connected sets of dimension $d$ are called \emph{cells} of the arrangement.
We refer the reader to \cite{EdeOroSei86,Ede87} for more information about hyperplane arrangements.


\subsubsection{The algorithm}
\label{sec algorithm}

\noindent
\textbf{Case $s = 1$.}
We set $\bar x^1 \in \bra{0,1}^{n_1}$ so that:
\begin{align}
\label{eq alg set last}
\bar x^1_k = 
\begin{cases}
1 & \text{if } \sum_{I \in \I} c^{I,1}_k > 0 \\
0 & \text{if } \sum_{I \in \I} c^{I,1}_k < 0
\end{cases}
\qquad \forall k \in [n_1],
\end{align}
meaning that in the case $\sum_{I \in \I} c^{I,1}_k = 0$, we can set arbitrarily $\bar x^1_k$ to either $0$ or $1$.
Return the solution $(\bar x^1)$.

\bigskip
\noindent
\textbf{Case $s \ge 2$.}
We partition the family $\I$ into $\I_\alpha$, $\I_\beta$, and $\I_\gamma$ as follows:
$\I_\alpha$ contains the sets $I \in \I$ that do not contain index $s$;
$\I_\beta$ contains the sets $I \in \I$ that strictly contain $\bra{s}$;
$\I_\gamma$ contains the sets $I \in \I$ that coincide with $\bra{s}$.
Note that $\abs{\I_\beta} = m_s$.
For every $I \in \I_\beta$, define variable $\lambda_I$.
For every $k \in [n_s]$, define the affine function
\begin{align}
\label{eq def h}
\begin{split}
h_k : & \quad \R^{\I_\beta} \to \R \\
& \quad \lambda \mapsto \sum_{I \in \I_\beta} \lambda_I c^{I,s}_k + \sum_{I \in \I_\gamma} c^{I,s}_k,
\end{split}
\end{align}
and the corresponding hyperplane in $\R^{\I_\beta}$:
\begin{align}
\label{eq hyperplanes}
H_k :=\bra{\lambda \in \R^{\I_\beta} : h_k(\lambda) = 0} \qquad \forall k \in [n_s].
\end{align}
Construct the arrangement of these hyperplanes in $\R^{\I_\beta}$ with the algorithm in \cite{EdeOroSei86,Ede87}, and denote by $\A$ the set of cells of the arrangement.

The remainder of the algorithm should be applied, separately, to each cell $C \in \A$, so we now fix one cell $C \in \A$.
The cell $C$ induces a signing of the hyperplanes $H_k$, for every $k \in [n_s]$.
Namely, for every $k \in [n_s]$, we know which inequality among
$h_k(\lambda) \ge 0$ and $h_k(\lambda) \le 0$ is valid for $C$, and note that precisely one of the two inequalities is valid for $C$, since $C$ is full-dimensional.
We then define the \emph{partial solution} corresponding to cell $C$ as $\bar x^s \in \bra{0,1}^{n_s}$ as follows:
\begin{align}
\label{eq alg set}
\bar x^s_k := 
\begin{cases}
1 & \text{if } h_k(\lambda) \ge 0 \text{ is valid for $C$} \\
0 & \text{if } h_k(\lambda) \le 0 \text{ is valid for $C$} 
\end{cases}
\qquad \forall k \in [n_s].
\end{align}
Define the \emph{child instance} of \cref{prob F} corresponding to cell $C$ as the optimization problem obtained, from \cref{prob F}, by fixing $x^s$ to $\bar x^s$:
\begin{align}
\label[problem]{prob F child}
\tag{F$_C$}
\begin{split}
\sum_{I \in \I_\gamma} \inner{c^{I,s}}{\bar x^s} + 
\max_{x^1, \dots, x^{s-1}} & \quad \sum_{I \in \I_\alpha} \prod_{j \in I} \inner{c^{I,j}}{x^j}
+ \sum_{I \in \I_\beta} \inner{c^{I,s}}{\bar x^s} \prod_{j \in I \setminus \bra{s}} \inner{c^{I,j}}{x^j}
\\
\st & \quad x^j \in \bra{0,1}^{n_j} \qquad \forall j \in [s-1].
\end{split}
\end{align}
\cref{prob F child} is essentially a new instance of \cref{prob F} with the parameter $s$ decreased by one.
However, there are two minor differences that we should point out.
First, there is a constant term $\sum_{I \in \I_\gamma} \inner{c^{I,s}}{\bar x^s}$ in the objective function, that we brought outside of the $\max$.
Second, 
there is a constant factor $\inner{c^{I,s}}{\bar x^s}$ in the sum over $I \in \I_\beta$.
This second issue can be easily remedied by absorbing the scalar $\inner{c^{I,s}}{\bar x^s}$ in precisely one of the subsequent inner products $\inner{c^{I,j}}{x^j}$, say the one corresponding to the smallest $j \in I \setminus \bra{s}$.
If we denote this index by $\hat j$, this is accomplished by redefining $c^{I,\hat j}$ as  $\inner{c^{I,s}}{\bar x^s} c^{I,\hat j}$ in the new instance.

We then apply the algorithm recursively to \cref{prob F child}.
Let $(\tilde x^1, \tilde x^2, \dots, \tilde x^{s-1})$ be the solution returned by the algorithm applied to \cref{prob F child}.
The algorithm then returns the solution $(\tilde x^1, \tilde x^2, \dots, \tilde x^{s-1}, \bar x^s)$ to \cref{prob F}.

\subsubsection{Correctness}

In this section, we show that the algorithm presented in \cref{sec algorithm} is correct, that is, that it returns an optimal solution to \cref{prob F}.

The proof is by induction on $s$.
In the base case we assume $s=1$, and so all sets $I \in \I$ are equal to $\bra{1}$.
\cref{prob F} takes the form
\begin{align*}
\max_{x^1} & \quad \sum_{I \in \I} \inner{c^{I,1}}{x^1} \\
\st & \quad x^1 \in \bra{0,1}^{n_1}.
\end{align*}
The objective function is linear in $x^1$ and can be written in the form 
\begin{align*}
\sum_{I \in \I} \inner{c^{I,1}}{x^1}
= \inner{\sum_{I \in \I} c^{I,1}}{x^1}
= \sum_{k \in [n_1]} \pare{\sum_{I \in \I} c^{I,1}_k} x^1_k.
\end{align*}
It is then simple to check that the solution $(\bar x^1)$ defined in \eqref{eq alg set last} and returned by the algorithm is optimal to \cref{prob F}.

Next, we consider the induction step, and we assume $s \ge 2$.
In our first claim, below, we write the objective function of \cref{prob F} as a linear function in $x^s$, if we think to all other variables $x^j$, for $j \in [s-1]$, as being fixed.
Let $\I_\alpha$, $\I_\beta$, and $\I_\gamma$ be the partition of the family $\I$ defined in the algorithm.

\begin{claim}
\label{claim obj}
The objective function of \cref{prob F} can be written in the form
\begin{align}
\label{eq obj}
\sum_{I \in \I_\alpha} \prod_{j \in I} \inner{c^{I,j}}{x^j}
+ \sum_{k \in [n_s]} h_k(\lambda) x^s_k,
\end{align}
where, for every $k \in [n_s]$, $h_k(\lambda)$ is the affine function defined by the algorithm in \eqref{eq def h}, and where
\begin{align*}
\lambda_I := \prod_{j \in I \setminus \bra{s}} \inner{c^{I,j}}{x^j} \qquad \forall I \in \I_\beta.
\end{align*}
\end{claim}

\begin{prf}
The objective function of \cref{prob F} can be written in the form
\begin{align}
\label{eq summands}
\sum_{I \in \I} \prod_{j \in I} \inner{c^{I,j}}{x^j}
= \sum_{I \in \I_\alpha} \prod_{j \in I} \inner{c^{I,j}}{x^j}
+ \sum_{I \in \I_\beta} \prod_{j \in I} \inner{c^{I,j}}{x^j}
+ \sum_{I \in \I_\gamma} \prod_{j \in I} \inner{c^{I,j}}{x^j}.
\end{align}
Since each $I \in \I_\beta$ strictly contains $\bra{s}$, we can write the second summand in \eqref{eq summands} as follows:
\begin{align*}
\sum_{I \in \I_\beta} \prod_{j \in I} \inner{c^{I,j}}{x^j}
& = \sum_{I \in \I_\beta} \pare{\prod_{j \in I \setminus \bra{s}} \inner{c^{I,j}}{x^j}} \inner{c^{I,s}}{x^s} \\
& = \inner{\sum_{I \in \I_\beta} \pare{\prod_{j \in I \setminus \bra{s}} \inner{c^{I,j}}{x^j}}c^{I,s}}{x^s} \\
& = \sum_{k \in [n_s]} \pare{\sum_{I \in \I_\beta} \prod_{j \in I \setminus \bra{s}} \inner{c^{I,j}}{x^j}c^{I,s}_k}{x^s_k}.
\end{align*}
Since each $I \in \I_\gamma$ coincides with $\bra{s}$, the third summand in \eqref{eq summands} can be written as follows:
\begin{align*}
\sum_{I \in \I_\gamma} \prod_{j \in I} \inner{c^{I,j}}{x^j}
& = \sum_{I \in \I_\gamma} \inner{c^{I,s}}{x^s} \\
& = \inner{\sum_{I \in \I_\gamma} c^{I,s}}{x^s} \\
& = \sum_{k \in [n_s]} \pare{\sum_{I \in \I_\gamma} c^{I,s}_k} x^s_k.
\end{align*}
From \eqref{eq summands}, we can then write the objective function of \cref{prob F} in the form
\begin{align*}
& \sum_{I \in \I_\alpha} \prod_{j \in I} \inner{c^{I,j}}{x^j}
+ \sum_{k \in [n_s]} \pare{\sum_{I \in \I_\beta} \prod_{j \in I \setminus \bra{s}} \inner{c^{I,j}}{x^j}c^{I,s}_k}{x^s_k}
+ \sum_{k \in [n_s]} \pare{\sum_{I \in \I_\gamma} c^{I,s}_k} x^s_k = \\
& = \sum_{I \in \I_\alpha} \prod_{j \in I} \inner{c^{I,j}}{x^j}
+ \sum_{k \in [n_s]} \sbra{\pare{\sum_{I \in \I_\beta} \prod_{j \in I \setminus \bra{s}} \inner{c^{I,j}}{x^j}c^{I,s}_k}+\pare{\sum_{I \in \I_\gamma} c^{I,s}_k}}{x^s_k} \\
& = \sum_{I \in \I_\alpha} \prod_{j \in I} \inner{c^{I,j}}{x^j}
+ \sum_{k \in [n_s]} \sbra{\sum_{I \in \I_\beta} \lambda_I c^{I,s}_k+\sum_{I \in \I_\gamma} c^{I,s}_k}{x^s_k} \\
& = \sum_{I \in \I_\alpha} \prod_{j \in I} \inner{c^{I,j}}{x^j}
+ \sum_{k \in [n_s]} h_k(\lambda) x^s_k.
\end{align*}
\end{prf}

Note that, in \eqref{eq obj}, the first sum and each $h_k(\lambda)$, for $k \in [n_s]$, generally depend on $x^j$, for $j \in [s-1]$, but do not depend on $x^s$.
\cref{claim obj} allows us to characterize the vector $x^s$ of all optimal solutions to \cref{prob F}, as we see in the next claim.

\begin{claim}
\label{claim property}
Let $\breve x$ be an optimal solution to \cref{prob F} and let 
\begin{align*}
\breve \lambda_I := \prod_{j \in I \setminus \bra{s}} \inner{c^{I,j}}{\breve x^j} \qquad \forall I \in \I_\beta.
\end{align*}
We have
\begin{align}
\label{eq claim conditions}
\breve x^s_k = 
\begin{cases}
1 & \text{if } h_k(\breve \lambda) > 0 \\
0 & \text{if } h_k(\breve \lambda) < 0.
\end{cases}
\qquad \forall k \in [n_s].
\end{align}
Furthermore, each solution obtained from $\breve x$ by flipping the value of any number of variables $\breve x^s_k$, $k \in [n_s]$, such that $h_k(\breve \lambda) = 0$, is optimal to \cref{prob F} as well.
\end{claim}

\begin{prf}
For a contradiction, assume there is $k \in [n_s]$ such that $\breve x^s_k$ does not satisfy the corresponding condition in \eqref{eq claim conditions}.
Consider now the solution obtained from $\breve x$ by flipping the value of variable $\breve x^s_k$, from $0$ to $1$, or from $1$ to $0$.
It follows from \cref{claim obj} that the objective value of this new solution is strictly larger than the objective value of $\breve x$.
This contradicts the optimality of $\breve x$, thus $\breve x$ must satisfy conditions \eqref{eq claim conditions}.

Next, let $k \in [n_s]$ such that $h_k(\breve \lambda) = 0$.
Consider the solution obtained from $\breve x$ by flipping the value of variable $\breve x^s_k$.
It follows from \cref{claim obj} that this new solution has the same objective value as $\breve x$, so it is optimal too.
\end{prf}


While there are $2^{n_s}$ possible binary vectors $x^s \in \bra{0,1}^{n_s}$, our goal is to exploit \cref{claim property} to identify a polynomial number of possible candidates.
The idea is to use hyperplane arrangements to consider all possible vectors $x^s$ that are compatible with \eqref{eq claim conditions}.
However, the function $h_k(\lambda)$, for every $k \in [n_s]$ depends on $x^1,x^2,\dots,x^{s-1}$ which amounts to $n_1+n_2+ \dots +n_{s-1}$ variables, and this would result in an exponential bound using hyperplane arrangements.
The key is to observe that $h_k(\lambda)$, in fact, depends only on the $\lambda_I$, for $I \in \I_\beta$, as suggested by our purposely chosen notation.

Consider the hyperplanes $H_k$, for every $k \in [n_s]$, defined by the algorithm in \eqref{eq hyperplanes}, and let $\A$ be the set of cells of the arrangement of these hyperplanes, 
as in the algorithm.
In the next claim, we show that there is at least one cell in $\A$ which provides a partial solution from an actual optimal solution to \cref{prob F}.

\begin{claim}
\label{claim match}
There exists an optimal solution $\mathring x$ to \cref{prob F} and a cell $C \in \A$ such that, if we denote by $\bar x^s \in \bra{0,1}^{n_s}$ the partial solution corresponding to cell $C$ defined in \eqref{eq alg set}, we have $\mathring x^s = \bar x^s$.
\end{claim}

\begin{prf}
Let $\breve x$ be an optimal solution to \cref{prob F} and let $\breve \lambda_I$, for every $I \in \I_\beta$, be defined as in \cref{claim property}.
Note that $\breve \lambda$ is in at least one cell in $\A$, and it may be contained in more than one.
Let $C$ be any cell in $\A$ containing $\breve \lambda$,
%
and let $\bar x^s \in \bra{0,1}^{n_s}$ be the partial solution corresponding to $C$ defined in \eqref{eq alg set}.
Now we compare $\breve x^s$ with $\bar x^s$ and show
\begin{align}
\label{eq compare}
\breve x^s_k = \bar x^s_k
\qquad \text{if } h_k(\breve \lambda) \neq 0 \qquad \forall k \in [n_s].
\end{align}
To prove \eqref{eq compare}, consider separately the cases $h_k(\breve \lambda) > 0$ and $h_k(\breve \lambda) < 0$.
In the first case we have $\breve x^s_k = 1$ from \eqref{eq claim conditions}.
Since $\breve \lambda \in C$, then $h_k(\lambda) \ge 0$ is valid for $C$, therefore we have $\bar x^s_k = 1$ from \eqref{eq alg set}.
In the second case we have $\breve x^s_k = 0$ from \eqref{eq claim conditions}.
Since $\breve \lambda \in C$, then $h_k(\lambda) \le 0$ is valid for $C$, therefore we have $\bar x^s_k = 0$ from \eqref{eq alg set}.
This concludes the proof of \eqref{eq compare}.

Let $\mathring x$ be the solution to \cref{prob F} defined as follows:
\begin{align*}
& \mathring x^j := \breve x^j \qquad \forall j \in [s-1] \\
& \mathring x^s_k := 
\begin{cases}
\breve x^s_k & \text{if } h_k(\breve \lambda) \neq 0 \\
\bar x^s_k & \text{if } h_k(\breve \lambda) = 0
\end{cases}
\qquad \forall k \in [n_s].
\end{align*}
From \cref{claim property}, $\mathring x$ is also an optimal solution to \cref{prob F}.
Due to \eqref{eq compare} and the definition of $\mathring x$, we have $\mathring x^s = \bar x^s$.
\end{prf}

Consider now the child instance \cref{prob F child} of \cref{prob F} corresponding to a cell $C$ from \cref{claim match}.
The next simple claim will be useful to conclude the proof of the correctness of our algorithm.

\begin{claim}
\label{claim simple}
Let $x^j \in \bra{0,1}^{n_j}$ for every $j \in [s-1]$.
The objective value of $(x^1, x^2, \dots, x^{s-1})$ in \cref{prob F child} equals the objective value of $(x^1, x^2, \dots, x^{s-1}, \bar x^s)$ in \cref{prob F}.
\end{claim}

\begin{prf}
Follows directly from the fact that \cref{prob F child} is obtained from \cref{prob F} by fixing $x^s$ to $\bar x^s$.
\end{prf}

Let $(\tilde x^1, \tilde x^2, \dots, \tilde x^{s-1})$ be the solution returned by the algorithm applied to \cref{prob F child}.

\begin{claim}
The solution $(\tilde x^1, \tilde x^2, \dots, \tilde x^{s-1}, \bar x^s)$ is optimal to \cref{prob F}.
\end{claim}

\begin{prf}
By induction, the solution $(\tilde x^1, \tilde x^2, \dots, \tilde x^{s-1})$ is optimal to \cref{prob F child}.
Let $\mathring x$ be an optimal solution to \cref{prob F} as in \cref{claim match}.
Since $(\tilde x^1, \tilde x^2, \dots, \tilde x^{s-1})$ is optimal to \cref{prob F child}, its objective value in \cref{prob F child} is greater than or equal to the objective value of $(\mathring x^1, \mathring x^2, \dots, \mathring x^{s-1})$ in \cref{prob F child}.
From \cref{claim simple}, the objective value of $(\tilde x^1, \tilde x^2, \dots, \tilde x^{s-1})$ in \cref{prob F child} equals the objective value of $(\tilde x^1, \tilde x^2, \dots, \tilde x^{s-1}, \bar x^s)$ in \cref{prob F}.
\cref{claim simple} also implies that the objective value of $(\mathring x^1, \mathring x^2, \dots, \mathring x^{s-1})$ in \cref{prob F child} equals the objective value of $(\mathring x^1, \mathring x^2, \dots, \mathring x^{s-1}, \bar x^s)$ in \cref{prob F}.
Hence, in \cref{prob F}, the objective value of $(\tilde x^1, \tilde x^2, \dots, \tilde x^{s-1}, \bar x^s)$ is greater than or equal to the objective value of $(\mathring x^1, \mathring x^2, \dots, \mathring x^{s-1}, \bar x^s)$.
From \cref{claim match}, we have $(\mathring x^1, \mathring x^2, \dots, \mathring x^{s-1}, \bar x^s) = (\mathring x^1, \mathring x^2, \dots, \mathring x^{s-1}, \mathring x^s)$, hence this solution is optimal to \cref{prob F}.
We then obtain that also $(\tilde x^1, \tilde x^2, \dots, \tilde x^{s-1}, \bar x^s)$ is optimal to \cref{prob F}.
\end{prf}


\subsubsection{Running time}


In this section, we prove the bound on the running time of our algorithm stated in \cref{th main ext}.

\begin{claim}
\label{claim runtime poly}
The algorithm presented in \cref{sec algorithm} performs at most
\begin{align*}
s \theta^{s-1} n_2^{m_2} n_3^{m_3} \cdots n_s^{m_s} \poly(n_1,n_2,\dots,n_s,\abs{\I})
\end{align*}
arithmetic operations, where $\theta$ is a constant.
\end{claim}

\begin{prf}
In the first iteration, the algorithm constructs at most $\theta n_s^{m_s}$ child instances of the form of \cref{prob F child} in at most $\theta n_s^{m_s}$ arithmetic operations, where $\theta$ is a constant \cite{EdeOroSei86,Ede87}.
The parameter $s$ associated with each child instance is decreased by one, and
the dimensions of the vectors of variables in the new instances are unchanged: $n_1,n_2,\dots,n_{s-1}$.
It follows from the construction of \cref{prob F child} that also the values $m_j$, for $j \in [s-1]$, associated with the child instances, are unchanged.
The total number of problems of the form \cref{prob F} considered in the recursive algorithm, and the number of arithmetic operations to construct them, is then at most
\begin{align*}
1+ \theta n_s^{m_s} + \theta^2 n_s^{m_s} n_{s-1}^{m_{s-1}} + \cdots + \theta^{s-1} n_s^{m_s} n_{s-1}^{m_{s-1}} \cdots n_2^{m_2}
\le 
s \theta^{s-1} n_s^{m_s} n_{s-1}^{m_{s-1}} \cdots n_2^{m_2},
\end{align*}
where the inequality holds because we can assume, without loss of generality, that $\theta \ge 1$.
For each instance considered, the number of remaining arithmetic operations performed by the algorithm is 
\begin{align*}
\poly(n_1,n_2,\dots,n_s,\abs{\I}).
\end{align*}
\end{prf}


To conclude the proof of \cref{th main ext},
it suffices to show the following result.

\begin{claim}
\label{claim runtime strong poly}
The size of the numbers produced in the course of the execution of the algorithm presented in \cref{sec algorithm} is polynomial in the size of the input.
\end{claim}

\begin{prf}
In each iteration of the algorithm, the size of the numbers produced by the algorithm is polynomial in the size of the data of the instance of \cref{prob F} considered at the beginning of that iteration.
This includes the numbers produced by the algorithm in \cite{EdeOroSei86,Ede87} to construct the arrangement of the hyperplanes.


On the other hand, we show that the size of the data of each child instance \cref{prob F child} constructed throughout the execution of the algorithm is polynomial in the size of the data of the original instance of \cref{prob F}.
In fact, from the construction of \cref{prob F child}, each child instance is obtained from the original instance of \cref{prob F} by fixing each component of $x^t, x^{t+1}, \dots, x^s$, for some $t \in \bra{2,3,\dots, s}$, to zero or one.
As a result, each number appearing in a child instance is obtaining by summing or multiplying polynomially many numbers from the original instance.
\end{prf}

\subsection{Corollaries}

In this section, we discuss some consequences of \cref{th main com}, for different optimization problems.

\subsubsection{An extension of factorized binary polynomial optimization}
First, we remark that our algorithm can be used, with minor modifications, to prove the same result for the slightly more general optimization problem obtained from \cref{prob F} by replacing each linear function in the objective with an affine function:
\begin{align}
\label[problem]{prob F'}
\tag{F$'$}
\begin{split}
\max_{x^1, \dots, x^s} & \quad \sum_{I \in \I'} \prod_{j \in I} \pare{\inner{c^{I,j}}{x^j} + d^{I,j}} \\
\st & \quad x^j \in \bra{0,1}^{n_j} \qquad \forall j \in [s].
\end{split}
\end{align}
The data of the problem consists of positive integers $s$ and $n_j$, for $j \in [s]$, 
a nonempty family $\I'$ of nonempty subsets of $[s]$,
vectors $c^{I,j} \in \Q^{n_j}$, and scalars $d^{I,j} \in \Q$,
for $I \in \I'$, $j \in I$.
We decided to give our algorithm only for \cref{prob F}, to avoid overloading the notation, which is already quite heavy.
Furthermore, the tractability of \cref{prob F'} can be obtained directly from the tractability of \cref{prob F}, as we see next.


\begin{corollary}
\label{cor affine}
\cref{prob F'} can be solved in strongly polynomial time, for any fixed $s$ and $\abs{\I'}$.
\end{corollary}

\begin{prf}
For every $I \in \I'$, the corresponding product in the objective function, $\prod_{j \in I} \pare{\inner{c^{I,j}}{x^j} + d^{I,j}}$, can be expanded to a sum
\begin{align*}
\sum_{I' \in S(I)} \prod_{j \in I} \inner{{c'}^{I,j}}{x^j} + d^I,
\end{align*}
where $S(I)$ is a set of nonempty subsets of $I$, ${c'}^{I,j} \in \Q^{n_j}$ for $j \in I$, and $d^I \in \Q$.
Therefore, \cref{prob F'} can be written in the form of \cref{prob F} with the same parameter $s$ and with $\abs{\I} \le 2^s \abs{\I'}$.
The result then follows from \cref{th main com}.
\end{prf}


\subsubsection{Explicit binary polynomial optimization on an $s$-partite hypergraph}

Using \cref{obs rank to factorization full}, we obtain the following corollaries to \cref{th main com} for \cref{prob E} on an $s$-partite hypergraph.

\begin{corollary}
Consider \cref{prob E} on an $s$-partite hypergraph.
This problem can be solved in strongly polynomial time if $s$ is fixed, and all the tensors of the coefficients of the objective function are given through a factorization of the form \eqref{eq factorization}, with $r$ fixed.
\end{corollary}

\begin{corollary}
Consider \cref{prob E} on an $s$-partite hypergraph.
This problem can be solved in strongly polynomial time if $s$ is fixed, and all the tensors of the coefficients of the objective function have fixed rank and rank factorizations are given.
\end{corollary}

\subsubsection{Binary tensor factorization and rank-1 Boolean tensor factorization}

Next, we discuss consequences of \cref{th main com} for the binary tensor factorization problem and the rank-1 Boolean tensor factorization problem.
\cref{th main com,prop BTF F,cor r1BTF F-U} directly imply the following results.

\begin{corollary}
\label{cor BTF factorization}
The rank-$t$ binary tensor factorization problem on an $n_1 \times n_2 \times \cdots \times n_s$ tensor
\begin{align*}
A = \sum_{p \in [q]} \pare{a^{p,1} \otimes a^{p,2} \otimes \cdots \otimes a^{p,s}},
\end{align*}
where $a^{p,j} \in \Q^{n_j}$, for every $p \in [q]$, $j \in [s]$, 
can be solved in strongly polynomial time for any fixed $q,s,t$.
\end{corollary}

\begin{corollary}
\label{cor r1BTF factorization}
The rank-1 Boolean tensor factorization problem on an $n_1 \times n_2 \times \cdots \times n_s$ tensor
\begin{align*}
A = \sum_{p \in [q]} \pare{a^{p,1} \otimes a^{p,2} \otimes \cdots \otimes a^{p,s}},
\end{align*}
where $a^{p,j} \in \Q^{n_j}$, for every $p \in [q]$, $j \in [s]$, 
can be solved in strongly polynomial time for any fixed $q,s$.
\end{corollary}

When the input tensor $A$ has fixed rank and a rank factorization is given, from 
\cref{cor BTF factorization,cor r1BTF factorization}, we obtain \cref{cor BTF rank,cor r1BTF rank}, stated in \cref{sec intro}.
In the special case $s=2$, \cref{cor r1BTF rank} implies \cref{cor r1BMF rank}, since a rank factorization of a matrix can be computed in strongly polynomial time via Gaussian elimination \cite{Edm67}.

\subsubsection{The quadratic case}
Regarding quadratic objective functions, \cref{cor main quad} can be obtained easily from \cref{th main com,obs rank to factorization full}.
We give a proof that does not use \cref{obs rank to factorization full}, since the algebra in this special case is significantly simpler.

\begin{prfc}[of \cref{cor main quad}]
Consider a matrix $Q \in \Q^{m \times n}$.
We can compute in strongly polynomial time, via Gaussian elimination \cite{Edm67}, a rank factorization of $Q$, that is, matrices $A \in \Q^{m \times r}$, $B \in \Q^{n \times r}$ such that $Q = AB^\transp$, where $r$ is the rank of $Q$ \cite{StewartBook}.
%
If we denote column $i$ of $A$ by $a^i$, and column $i$ of $B$ by $b^i$, for every $i \in [r]$, we can write
\begin{align*}
x^\transp Q y 
& = x^\transp A B^\transp y \\
& =  x^\transp \pare{ \sum_{i \in [r]} \pare{a^i \otimes {b^i}} } y \\
& =  x^\transp \pare{ \sum_{i \in [r]} \pare{a^i {b^i}^\transp} } y \\
& = \sum_{i \in [r]} \pare{x^\transp a^i {b^i}^\transp y} \\
& = \sum_{i \in [r]} \pare{\inner{a^i}x \inner{b^i}y}.
\end{align*}

We apply the above argument to each matrix $Q^{i,j} \in \Q^{n_i \times n_j}$, for $i,j \in [s]$ with $i<j$, in the objective function of \cref{prob Q}.
There are $(s^2-s)/2$ of these matrices, and we denote by $r$ the maximum of their ranks.
\cref{prob Q} can then be written in the form of \cref{prob F}, where $\abs{\I}$ is at most $r(s^2-s)/2+s$.
This number is fixed, since by assumption $s$ and $r$ are fixed.
The result then follows from \cref{th main com}.
\end{prfc}

\section{Comparison with known tractability results}
\label{sec comparison}

Every instance of \cref{prob F} can be reformulated as an instance of \cref{prob E} by expanding all products in the objective function; thus we can employ known algorithms for \cref{prob E} to solve \cref{prob F}.
In this section we show that our \cref{th main com} does not follow from known tractability results of \cref{prob E}.
To the best of our knowledge, five main polynomially solvable classes of \cref{prob E} on a hypergraph $H=(V,E)$ have been identified so far.
These are instances such that:
\begin{enumerate}[noitemsep,nolistsep,label=C\arabic*,leftmargin=*]
\item
\label{class 1}
The objective function is supermodular (see Chapter 45 in \cite{SchBookCO});
\item
\label{class 2}
$H$ is a $\beta$-acyclic hypergraph \cite{dPDiG22SODA,dPDiG23ALG,dPKha23MPA} (see also \cite{dPKha18SIOPT,BucCraRod18,dPKha21MOR});
\item
\label{class 3}
$H$ is a cycle hypergraph \cite{dPDiG21IJO}; 
\item
\label{class 4}
The primal treewidth of $H$ is bounded by $\log(\poly(\abs{V},\abs{E}))$ 
\cite{CraHanJau90,Lau09,BieMun18};
\item
\label{class 5}
The incidence treewidth of $H$ is bounded by $\log(\poly(\abs{V},\abs{E}))$ 
\cite{CapdPDiG24}.
\end{enumerate}

It will suffice to consider the special case of \cref{prob F} in \cref{ex cu}: \cref{prob F-U} with $r=1$, and with each vector $c^{1,j} \in \Q^{n_j}$, for every $j \in [s]$, having all nonzero components.
It is also sufficient to consider the case $n_1=n_2=\cdots=n_s$, and we set $n:=n_1$.
\cref{th main com} implies that this problem can be solved in strongly polynomial time for any fixed $s$.
As we saw in \cref{ex cu}, expanding all products in the objective function, the problem takes the form of \cref{prob E-U}, where the corresponding hypergraph is $s$-partite, with sides $V^j$, for $j \in [s]$, of cardinality $n$, is $s$-uniform, and contains all possible $n^s$ edges;
We denote this hypergraph by $H^s_n$.

The following facts imply that 
even the polynomial solvability of this special case of \cref{prob F}, for every $s \ge 2$ fixed, does not follow from any of the known tractable classes \ref{class 1}--\ref{class 5} above.

\begin{fact}
\label{fact 1}
The objective functions of the problems considered in \cref{ex cu} are not generally supermodular, for every $s \ge 2$, $n \ge 1$.
\end{fact}

\begin{fact}
\label{fact 2}
The hypergraph $H^s_n$ is not $\beta$-acyclic, for every $s \ge 2$, $n \ge 2$.
\end{fact}

\begin{fact}
\label{fact 3}
The hypergraph $H^s_n$ is not a cycle hypergraph, for every $s \ge 2$, $n \ge 3$.
\end{fact}

\begin{fact}
\label{fact 4}
The primal treewidth of $H^s_n$ is $(s-1) n$, for every $s \ge 2$, $n \ge 1$.
\end{fact}

\begin{fact}
\label{fact 5}
The incidence treewidth of $H^s_n$ is at least $n$, for every $s \ge 2$, $n \ge 1$.
\end{fact}

\cref{fact 1,fact 2,fact 3} are easy to verify.

\begin{prfc}[of \cref{fact 1}]
It suffices to partition the set $[s]$ into nonempty sets $S_1, S_2$, and construct one solution $\bar x = (\bar x^1, \bar x^2, \dots, \bar x^s)$ with all $\bar x^j$, $j \in S^1$, with all components one, and all $x^j$, $j \in S^2$, with all components zero.
The objective values of $\bar x$ and $1-\bar x$ are zero, and so is the objective value of the componentwise minimum $\bar x \downarrow (1-\bar x)$.
On the other hand, the objective value of the componentwise maximum $\bar x \uparrow (1-\bar x)$ is the sum of all edge costs, which can easily be made negative.
\end{prfc}

\begin{prfc}[of \cref{fact 2}]
Let $v^j_1,v^j_2$ be two distinct nodes in $V^j$, for every $j \in [s]$.
A $\beta$-cycle of length four is 
\begin{align*}
v^1_1, 
\bra{v^1_1,v^2_1,\dots,v^s_1},
v^s_1,
\bra{v^1_2,v^2_1,v^3_1,\dots,v^s_1}, 
v^1_2,
\bra{v^1_2,v^2_2,\dots,v^s_2}, 
v^s_2,
\bra{v^1_1,v^2_2,v^3_2,\dots,v^s_2},
v^1_1.
\end{align*}
\end{prfc}

\begin{prfc}[of \cref{fact 3}]
For every $s \ge 2$ and $n \ge 3$, in the hypergraph $H^s_n$, each node is contained in at least three edges.
However, in a cycle hypergraph, each node is contained in at most two edges.
\end{prfc}

In the remainder of the section, we prove \cref{fact 4,fact 5}.
We start by defining the primal treewidth and the incidence treewidth of a hypergraph $H=(V,E)$.
First, we associate two graphs to $H$. 
The \emph{primal graph} $G_{\prim}(H)$ of $H$ is defined as the graph whose node set is $V$ and edge set is $\{\{u,v\} : u \neq v, \ \exists e \in E, \ \{u,v\} \in e \}$. 
Intuitively, the primal graph is obtained by replacing every edge of $H$ by a clique. 
The \emph{incidence graph} $G_{\inc}(H)$ of $H$ is defined as the bipartite graph whose node set is $V \cup E$ and the edge set is $\{\{v,e\} : v \in V, \ e \in E, \ v \in e\}$. 
The \emph{primal treewidth} $\ptw(H)$ of $H$ is the treewidth of its primal graph, that is $\ptw(H) = \tw(G_{\prim}(H))$, while the \emph{incidence treewidth} $\itw(H)$ is the treewidth of its incidence graph, that is $\itw(H)=\tw(G_{\inc}(H))$. 
In the next result, we characterize the treewidth of a \emph{complete $s$-partite graph,} which is an $s$-partite graph that contains all possible edges (of cardinality two).

\begin{lemma}
\label{lem tw complete}
Let $K$ be the complete $s$-partite graph, with $s \ge 2$, and with sides $V^j$ of cardinality $n$, for $j \in [s]$.
Then $\tw(K) = (s-1) n$.
\end{lemma}

\begin{prf}
It is well-known that the treewidth of $K$ is the minimum size of the largest clique minus one, in a chordal completion of $K$.

Let $\ind \in [s]$.
Let $G$ be the graph obtained from $K$ by adding, for each $j \in [s] \setminus \bra{\ind}$, all edges between all pairs of nodes in $V^j$.
The graph $G$ is chordal because it has the following perfect elimination ordering: First we list all nodes in $V^{\ind}$ in any order, and next all nodes in $V \setminus V^{\ind}$ in any order.
Since any maximal clique of $G$ consists of all nodes in $V \setminus V^{\ind}$ and one node in $V^{\ind}$, the size of the largest clique in $G$ is $(s-1)n + 1$.

Now let $G'=(V,F)$ be a chordal graph containing $K$.
We claim that there is $\ind \in [s]$ such that, for every $j \in [s] \setminus \bra{\ind}$, $F$ contains all edges between all pairs of nodes in $V^j$.
Assume by contradiction that this does not hold.
Then there exist two distinct indices $s,t \in [s]$, two nodes $u_1,u_2 \in V^s$ and two nodes $v_1,v_2 \in V^t$ such that $\{u_1, u_2\}, \{v_1, v_2\} \notin F$.
But then $u_1,v_1,u_2,v_2,u_1$ is a chordless cycle of $G'$ of length 4, a contradiction. 
This completes the proof of our claim.
Hence, $G'$ contains as a subgraph the graph $G$ from the previous paragraph.
So the size of the largest clique in $G'$ is at least $(s-1)n + 1$.
The treewidth of $K$ is then $(s-1)n$.
\end{prf}

We are now ready to prove \cref{fact 4,fact 5}.
\begin{prfc}[of \cref{fact 4}]
The graph $G_{\prim}(H^s_n)$ is the complete $s$-partite graph with sides $V^j$ of cardinality $n$, for $j \in [s]$.
From \cref{lem tw complete}, $\tw(G_{\prim}(H)) = (s-1)n$.
We obtain $\ptw(H) = \tw(G_{\prim}(H))=(s-1)n$.
\end{prfc}

\begin{prfc}[of \cref{fact 5}]
Consider the graph $G_{\inc}(H^s_n)$, whose node set is $V \cup E$.
Let $G'$ be the minor of $G_{\inc}(H^s_n)$ obtained by deleting all nodes in $V^3 \cup V^4 \cup \cdots \cup V^s$, and then contracting the edges with one node in $E$ and the other node in $V^2$.
$G'$ is then the complete $2$-partite graph with sides $V^1$, $V^2$.
From \cref{lem tw complete}, $\tw(G') = n$.
We obtain $\itw(H^s_n)=\tw(G_{\inc}(H^s_n)) \ge \tw(G') = n$,
where the inequality holds because $G'$ is a minor of $G_{\inc}(H^s_n)$ (see, e.g., lemma~14 in \cite{BolKos11}).
\end{prfc}


\paragraph{The quadratic case.}


The special case $s=2$ of \cref{prob Q} is the \emph{binary bipartite quadratic optimization} (BQO) problem studied in \cite{PunSriKar15a,PunSriKar15b,Pun22}:

\begin{align}
\label[problem]{prob BQO}
\tag{BQO}
\begin{split}
\max_{x^1, x^2} & \quad {x^1}^\transp Q {x^2} + {c^1}^\transp x^1 + {c^2}^\transp x^2 \\
\st & \quad x^1 \in \bra{0,1}^{n_1}, \ x^2 \in \bra{0,1}^{n_2},
\end{split}
\end{align}
where $Q \in \Q^{n_1 \times n_2}$, $c^1 \in \Q^{n_1}$, and $c^2 \in \Q^{n_2}$.
%
%
A polynomial time algorithm for \cref{prob BQO}, under the assumption that the rank of $Q$ is fixed, follows by combining the algorithm for the continuous relaxation of \cref{prob BQO} in \cite{HlaCerRad21} with the rounding procedure in \cite{PunSriKar15b}, as observed in \cite{Pun22}.
This algorithm is not strongly polynomial, since it needs to solve linear optimization problems.
An algorithm that solves \cref{prob BQO} in strongly polynomial time if $Q$ has rank one, is given in \cite{PunSriKar15a,Pun22}.
In \cite{PunSriKar15a}, the authors also present an algorithm for \cref{prob BQO} that is strongly polynomial if the rank of $Q$ is fixed and some ``dual non-degeneracy assumptions'' are satisfied.
As mentioned by the authors, it seems possible to lift these dual non-degeneracy assumptions by constructing an appropriate perturbation of the objective function.

The special case $s=2$ of our \cref{cor main quad}, implies that \cref{prob BQO} can be solved in strongly polynomial time if the rank of $Q$ is fixed.
Therefore, our result, even in the very special setting, significantly expands previously known results, both in terms of generality, since our graph is $s$-partite and not just bipartite, and in terms of computational complexity, since our algorithm is strongly polynomial rather than just weakly polynomial.


\bigskip

\begin{small}
\noindent
\textbf{Funding: }Alberto~Del~Pia is partially funded by AFOSR grant FA9550-23-1-0433. Any opinions, findings, and conclusions or recommendations expressed in this material are those of the authors and do not necessarily reflect the views of the Air Force Office of Scientific Research.
\end{small}

\ifthenelse {\boolean{MPA}}
{
\bibliographystyle{spmpsci}
}
{
\bibliographystyle{plain}
}


\end{document}